
\magnification=\magstep1
\baselineskip=20truept
\parindent=0pt
\null\vskip1truein
\font\tw=cmr12
\def\pne{{\cal P}({}^nE)}

\def\png{{\cal P}({}^nG)}
\def\pme{{\cal P}({}^mE)}
\def\C{\bf C}
\def\o{\over}
\def\dis{\displaystyle}
\def\sumi#1#2{\sum_{i=#1}^{#2}}
\def\sumj#1#2{\sum_{j=#1}^{#2}}
\def\suml#1#2{\sum_{l=#1}^{#2}}
\def\summ#1#2{\sum_{m=#1}^{#2}}
\def\sumn#1#2{\sum_{n=#1}^{#2}}
\def\onpis{\widehat{\bigotimes
_{n,s,\pi}}}
%
%
\centerline{{\tw $Q$-Reflexive Banach Spaces}}
\centerline{by}
\centerline{\tw Richard M. Aron\footnote*{{\rm The research
for this article was undertaken during the academic year
1992--1993 while R.~M.~Aron was visiting professor at
University College Dublin. Research supported in part by
NSF Grant Int-9023951.}} and Se\'an Dineen}
\centerline{\tw (Kent State University)\qquad (University
College Dublin)}
\bigskip\bigskip
Let $E$ be a Banach space. There are several natural ways in
which any polynomial $P\in\pne$ can be extended to $\tilde
P\in {\cal P}({}^nE'')$, in such a way that the extension
mapping is continuous and linear (see, for example, [6]).
Taking the double transpose of the extension mapping $P \to
\tilde P$ yields a linear, continuous mapping from $\pne''$
into ${\cal P}({}^nE'')''$.  Further, since ${\cal
P}({}^nE'')$ is a dual space, it follows that there is a
natural projection of ${\cal P}({}^nE'')''$ onto ${\cal
P}({}^nE'')$, and thus we have a  mapping of $\pne''$ into
${\cal P}({}^nE'')$. If all polynomials on a Banach space
$E$ are weakly continuous on bounded sets then these
mappings from $\pne''$ into ${\cal P}({}^nE'')$ all
coincide, and have a particularly simple description. We
discuss this in some detail below.

In this article we restrict ourselves to the situation in
which all polynomials on $E$ are weakly continuous on
bounded sets, and we study when this mapping is an
isomorphism. As we will see, if three ``ingredients'' are
present then the mapping will be an isomorphism: (1) $E''$
has the Radon-Nikodym property [16], (2) $E''$ has the
approximation property [27], and (3) every polynomial on $E$
is weakly continuous on bounded sets. In addition, we will
construct an example of a quasi-reflexive (non-reflexive)
Banach space $E$ for which the extension mapping is an
isomorphism.

It is well known that $\pne$ is isomorphic to
$(\bigotimes\limits_{n,s}E)'$, the dual of the $n$-fold {\it
symmetric} tensor product of $E$ endowed with the projective
topology. In fact, our results carry over to the space
$(\bigotimes\limits_{n}E)'$. However, since our interest is
in polynomials and holomorphic functions on $E$, we have
preferred to concentrate on symmetric tensor products.

By [9,theorem 2.9] if all continuous polynomials are weakly
continuous on bounded sets then they are in fact uniformly
weakly continuous on bounded sets and so have a unique
extension to polynomials on $E''$ which are weak$^*$
continuous on bounded sets. The mapping $P\in\pne\to\tilde
P(x)$, where $x\in E''$ and $\tilde P$ is the weak$^*$
continuous extension of $P$, is continuous when spaces of
polynomials are given their norm topology.  We define the
canonical mapping of the symmetric tensor product of $E''$
into the dual of $\pne$ in the following fashion
$$\eqalign{
J_n:\bigotimes\limits_{n,s}E''&\longrightarrow\pne'\cr
\bigotimes\limits_{n}x&\longrightarrow\bigl
[P\in\pne\to\tilde
P(x)\bigr]\cr}$$
Since $\|P\|=\|\tilde P\|$ it follows that
$$\eqalign{
\|J_n(\bigotimes_n x)\bigr\|&=\sup_{\scriptstyle P\in\pne
\atop\scriptstyle \|P\|\le1}\bigl|\tilde P(x)\bigr|\le
\sup_{\scriptstyle P\in\pne
\atop\scriptstyle \|P\|\le1}\|\tilde P\|\cdot\|x\|^n\cr
&\le\|x\|^n\cr}$$
Hence
$$\eqalign{
\bigl\|J_n(z)\bigr\|&\le\inf\left\{\sumi 1 m\|x_i\|^n;z
=\sumi 1
m\bigotimes_nx_i\right\}\cr
&=\|z\|_\pi\cr}$$
and $J_n$ can be extended in a unique fashion
to a continuous
linear mapping
$$J_n:\onpis E''\to\pne'$$
The transpose ${}^tJ_n$ is the required canonical mapping
 from
$$\pne''\quad\hbox{ into }{\cal P}
({}^nE'')\simeq\left(\bigotimes%
_{n,s,\pi}E''\right)'$$
Various attempts at defining a canonical mapping from
$\pne''$ into ${\cal P}({}^nE'')$ have convinced us that the
class of Banach spaces we are considering is the natural
class \footnote{(1)}{M.~Gonzalez has recently obtained a
characterization of $Q$-reflexive spaces which justifies our
choice.} for such a mapping.  Note that if $J$ denotes the
canonical mapping from $E$ into $E''$ then $J_1\circ J=J$.

\proclaim Definition 1. We shall say that a Banach space $E$
is $Q$-reflexive if the canonical mapping ${}^tJ_n$ is an
isomorphism from $\pne''$ onto ${\cal P}({}^nE'')$ for all
$n$.

It is easily seen that a reflexive Banach space is
$Q$-reflexive if and only if $\pne$ is reflexive for all
$n$.  This collection of spaces has been studied in
[1,2,3,4,18,19,21,23] and so we confine ourselves here to
nonreflexive Banach spaces which are $Q$-reflexive.  If $E$
is $Q$-reflexive then, since $\pne$ is a dual space it
follows that $\pne$ is $1$-complemented in ${\cal
P}({}^nE'')$ and we have the decomposition
$${\cal P}({}^nE'')={\cal P}_{\omega*}({}^nE'')\oplus {\cal
P}_{E^\perp}({}^nE'')$$ where ${\cal
P}_{\omega*}({}^nE'')\approx {\cal P}({}^nE)$ denotes the
set of all polynomials on $E''$ which are weak* continuous
on bounded sets and
$${\cal P}_{E^\perp}({}^nE'')=\left\{P\in {\cal
P}({}^nE'');P|_E\equiv0\right\}.$$

The following result of Gutierrez, which improves a result
in [7], reduces the study of $Q$-reflexive spaces to a more
manageable collection of Banach spaces in which we find a
more practical characterization of polynomials which are
weakly continuous on bounded sets.

\proclaim Proposition 2 ([23]).  If $E$ is a
Banach space which contains a copy of $l_1,$ then $E$ admits
a $\C$-valued homogeneous polynomial which is not weakly
continuous on bounded sets.

\proclaim Proposition 3.  If $E$ is a Banach space and
$l_1\not\hookrightarrow E$ then the following are
equivalent:\hfill\break
(a) all continuous polynomials on $E$ are weakly continuous
on bounded sets\hfill\break
(b) all continuous polynomials on $E$ are weakly
sequentially continuous at the origin.

\noindent{\bf Proof}.  By [9, proposition 2.12] the
continuous polynomials on~$E$ which are weakly sequentially
continuous are  weakly (uniformly) continuous on bounded
sets. By~[3] it is sufficient to check weak sequential
continuity at the origin.
\bigskip
To obtain  examples of non-reflexive $Q$-reflexive Banach
spaces we note that the proof of theorem~1 in~[18] does not
require reflexivity and consequently the following is true.

\proclaim Proposition 4.  If $E$ is a Banach space such that
no spreading model built on a normalised weakly null
sequence has a lower $q$-estimate for any $q<\infty$ then
any continuous polynomial on any subspace of $E$ is weakly
sequentially continuous at the origin.

The dual of a Banach space $E$ has the Radon-Nikodym
Property (RNP) if and only if each separable subspace of $E$
has a separable dual [16].  Such spaces are also called
Asplund space.  If $E$ is an Asplund space then
$l_1\not\hookrightarrow E$ since $(l_1)'=l_\infty$ is
nonseparable.  If $E'$ is an Asplund space and
$l_1\hookrightarrow E$ then by [15,p211],
$l_1\hookrightarrow l_1(2^N)\hookrightarrow E'$ and
consequently $E'$ contains a separable subspace with
nonseparable dual.  This is impossible and hence if either
$E$ or $E'$ is Asplund it follows that
$l_1\not\hookrightarrow E$.  

\bigskip

The authors are grateful to David Yost for pointing out an
error in the original version of the following result and
for informing them of reference [29]. 

\proclaim Theorem 5.  Let $E$ denote a Banach space such
that no spreading model built on a normalised weakly null
sequence has a lower $q$-estimate for any $q<\infty$.

{\leftskip=35truept\rightskip35truept
\sl
{(a)}  If $E'$ has RNP and the approximation property
then $\pne$ has RNP for all $n$.\hfill\break
{(b)}  If $E''$ has RNP and the approximation property
then $E$ is a $Q$-reflexive Banach space.}
\bigskip

{\bf Proof}.  By the above remarks it follows that in both
cases $l_1\not\hookrightarrow E$. By propositions~3 and~4 we
see that all continuous polynomials on any subspace
of $E$ are weakly
continuous on bounded sets.  By~[27, proposition~1.e.7] and
our hypothesis in (a) it follows that $E'$ has the
approximation property in both cases.  By [9,
corollary~2.11] it follows that $\dis\png\simeq\left(
\widehat{\bigotimes_{n,s,\epsilon}}G'\right)$, 
for any subspace $G$ of $E$ such that $G'$ has the
approximation property. We now
complete the proof of (a).

By [16, p. 218], it suffices to show that any separable
subspace $H$ of $\pne$ is isomorphic to a subspace of a separable dual
space. Suppose that $\{\phi_j^n\}_{j=1}^{\infty}, \ \	\phi_j 
\in E'$, spans a dense subspace of $H$. Let $F$ denote the
closed subspace of $E'$ generated by $\{\phi_j\}_{j=1}^{\infty}.$
By [29, Proposition 2], there exists a separable subspace $E_1$ of $E$ 
and a complemented subspace $F_1$ of $E'$ such that $F \subset
F_1$ and $E_1'\simeq F_1.$ Since  $E'$ has the RNP and the
approximation property and $F_1$ is complemented in $E'$, it
follows that $E_1' \simeq F_1$ also has both of these properties.
Hence $$H \subset   \left(
\widehat{\bigotimes_{n,s,\epsilon}}F_1\right) \simeq
\left(
\widehat{\bigotimes_{n,s,\epsilon}}E_1'\right) \simeq {\cal P}({}^nE_1) .$$
This implies that ${\cal P}({}^nE_1)$ is a separable dual space and that 
$H$ is isomorphic to a subspace of a separable dual space. Thus, $\pne$ 
has RNP. This completes the proof of (a). 




We now complete the proof of (b). By the above
$\dis\pne=\left(
\widehat{\bigotimes_{n,\epsilon,s}}E'\right)$. Since $E''$
has the RNP and the approximation property we have
by~[14,22] $$\pne'\simeq\left(\widehat{
\bigotimes_{n,\epsilon,s}}E'\right)'
\simeq\widehat{\bigotimes_{n,\pi,s}}E''$$ where the
isomorphism $I_n$ between these spaces satisfies
$$\bigl(I_n(\bigotimes_{n,s}x'')\bigr)(\phi^n)=
\bigl(x''(\phi)\bigr)^n$$ for all $x''\in E''$ and $\phi\in
E'$.

Hence $I_n=J_n$ and ${}^tJ_n$ is an isomorphism, i.e.~$E$ is
$Q$-reflexive.  This completes the proof.

\proclaim Corollary 6.  If $E$ is a Banach space such that
no spreading model built on a normalised weakly null
sequence has a lower $q$-estimate for any $q<\infty$ and
$E'$ has the approximation property then
$\dis\widehat{\bigotimes_{n,\pi,s}}E$ is Asplund for all $n$
if and only if $E$ is Asplund.

{\bf Example 7}.  Since all spreading models built on a
normalised weakly null sequence in $c_0$ are isomorphic to
$c_0$ $\bigl([\rm 9,p. 72]\bigr)$ and $c_0'=l_1$ has the
approximation property and RNP it follows that ${\cal
P}(^nc_0)$ has RNP for all $n$. This, however, also follows
immediately from the fact, proved in [5], that ${\cal
P}(^nc_0)$ is separable and thus as a separable dual space
${\cal P}(^nc_0)$ has RNP.  We now show that $c_0$ is not
$Q$-reflexive.  In this case the canonical mapping is
$J_n:{\cal P}_N(^nl_1)\longrightarrow{\cal P}_I(^nl_1)$
where ${\cal P}_N(^nl_1)$ and ${\cal P}_I(^nl_1)$ are
respectively the $n$-homogeneous nuclear and integral
polynomials on $l_1$.  There are various ways in which one
can show that $J_n$ is not an isomorphism and hence that
$c_0$ is not $Q$-reflexive.  For instance
$\bigl(\overline{B_{l_\infty}},
\sigma(l_\infty,l_1)\bigr)\simeq\Delta^N$ where $\Delta$ is
the closed unit disc in $\C$.  Let $\mu$ denote a Borel
probability measure on $\Delta$ such that $\dis\int_\Delta
z\,d\mu(z)=0$ and $\dis\int_\Delta z^2\,d\mu(z)=1$.  If
$\dis v=\prod_{n=1}^\infty\mu_n$ on $\Delta^N$, where
$\mu_n=\mu$ all $n$,  then the mapping
$$\eqalignno{
(x_n)_n\in
l_1&\to\int_{[0,1]^N}\bigl[(y_n)_n
\bigl((x_n)_n\bigr)\bigr]^2\,dv
\bigl((y_n)_n\bigr)\cr
&=\sum_{n,m}x_nx_m\int_{\Delta^N}
y_ny_m\,dv\bigl((y_n)_n\bigr)&(*)\cr
&=\sum_{n=1}^\infty x_n^2\cr}$$
defines an element of ${\cal P}(^nc_0)'={\cal
P}_I(^nl_1)$.  The associated integral mapping from $l_1$
into $l_\infty$ is not compact and hence not nuclear.  This
proves our claim.

{\bf Remarks}. (a) The action of the polynomial represented
by (*) on ${\cal P}(^2c_0)$ is given by
$$P\in{\cal P}(^2c_0)\to\sum_{n=1}^\infty P(e_n)$$
where $(e_n)_n$ is the standard unit vector basis in $c_0$.
Cleasly, replacing $\mu_n$ by a point mass at the origin
for all $x\not\in M$, $M$ some subset of $N$, we see that
$\sum_{n\in M}P(e_n)<\infty$ for all $P\in{\cal P}(^2c_0)$.
This provides another proof of a result in [8,13,31], namely
that $\sum_{n=1}^\infty\bigl|P(e_n)\bigr|<\infty$ for all
$P\in{\cal P}(^2(c_0)$.

(b) The above shows that there is a one to one
correspondence between ${\cal P}_I(^2l_1)$ and the
convariances of signed Borel measures on $\Delta^N$.

By using Grothendieck's inequality it follows that
$l_\infty\widehat{\bigotimes}_\pi l_\infty$ can be
represented as the set of all functions of the form
$$f:N\times N\to\C$$
where $f(n,m)=\bigl\langle g(n),h(m)\big\rangle$
and $\bigl(g(n)\bigr)_n$ and $\bigl(h(m)\bigr)_m$ are
relatively compact sequences  in $l_2$.

The space $l_1\widehat{\bigotimes}_\epsilon l_1$ has a
representation as the set of all series of the form
$\dis\sum_{n=1}^\infty e_n^*\bigotimes x_n$ where
$(e^*_n)_n$ is the standard unit vector basis in $l_1$ and
$(x_n)_n$ is an unconditionally convergent series in $l_1$
([2,20]).

(c)  Example 7 shows that the hypothesis ``$E''$ has RNP''
cannot be removed in Theorem~5 (b).

Our next step is to produce an example of a non-reflexive
$Q$-reflexive Banach space.  Theorem~5 clearly suggests that
we should consider a space like Tsirelson's space but not a
reflexive space and we were thus led to the
(quasi-reflexive) James space modeled on the original
Tsirelson space.  It is clear, however, that this is
representative of a class of examples which can be found
using the appropriate properties of the James and Tsirelson
space.  We refer to [27] for details concerning the James
spaces and the dual of the original Tsirelson space.
Following standard notation we let $T$ denote the dual of
the original Tsirelson space normed as in [27, p.~95] and we
denote its dual by $T^*$.  The James space modelled on $T$
is discussed in [11] and [26].

Let $(t_n)_{n=1}^\infty$ denote the standard unconditional
basis for $T^*$ which is dual to the basis given in [27,
p.~95].  The following two properties of $T^*$ play an
essential r\^ole in our construction.  For the sake of
completeness we include a proof of~(2).

For any positive integer $n$ we have
$$\left\|\sum_{j=n}^{2n}a_it_j\right\|_{T^*}\le2\sup_{n\le
i\le2n}|a_i|\eqno(1)$$
$\bigl($[12, proposition 1.7]$\bigr)$.

If $(k_j)_{j=1}^\infty$ is an increasing sequence of
integers, $k_1=0$, then for any
$\sum_{j=1}^\infty a_jt_j\in
T^*$ we have
$$\left\|\sum_{j=1}^\infty a_jt_j\right\|_{T^*}\le
\left\|\sum_{j=1}^\infty\left\|
\sum_{i=k_{j}+1}^{k_{j+1}}a_it_i
\right\|_{T^*}t_j\right\|_{T^*}\eqno(2)$$
([12, lemma II.1 and notes and remarks p22]).

\noindent{ Proof} of inequality (2).
Let
$$u_j={\sumi {{k_j+1}}
{{k_{j+1}}}a_it_i\o\left\|\sumi {{k_{j+1}}}
{{k_{j+1}}} a_it_i\right\|_{T^*}}.$$
Then
$$\eqalign{
\left\|\sum_{j=1}^\infty a_jt_j\right\|_{T{}^*}&=
\left\|\sumj 1 \infty\left(
\sum_{i=k_j+1}^{k_{j+1}}a_jt_j\right)
\right\|_{T{}^*}\cr
&=\left\|\sumj 1 \infty\left\|\sumi {{k_{j}+1}}
{{k_{j+1}}}a_it_i\right\|_{T^*}
u_j
\right\|_{T{}^*}\cr
&\le\left\|\sumj 1 \infty\left\|\sumi {{k_j+1}}
{{k_{j+1}}}a_it_i\right\|_{T^*}
t_{k_j+1}\right\|_{T^*}\cr
&\le\left\|\sumj 1 \infty\left\|\sumi {{k_{j}+1}}
{{k_{j+1}}}a_it_i\right\|_{T^*}
t_{j}\right\|_{T^*}\cr}$$
where we use $u_j$ to get the normalised sequence in [12]
and we have reversed the inequality since we are dealing
with $T{}^*$ in place of $T$. Since $j\le k_j+1$ and moving
the support to the right in~$T$ increases the norm ([12,
proposition 1.9(3)]) it follows that in $T{}^*$ moving the
support to the left increases the norm and this is what we
have done here.

For $(a_n)_{n=1}^\infty\in c_{00}$, the space of all
sequences which are eventually zero.  Let
$$\bigl\|(a_n)_{n=1}^\infty\bigr\|_{T{}_J^*}=
\sup_{\scriptstyle p_1<p_2<\cdots<p_{2k}
\atop\scriptstyle  k}
\left\|\sum_{j=1}^k(a_{p_{2j-1}}-
a_{p_{2j}})t_j\right\|_{T^*}\eqno(3)$$
The completion of $c_{00}$ with respect to the norm
$\|\>\cdot\>\|_{T{}_J^*}$ is denoted by ${T{}_J^*}$ and is
called the Tsirelson${}^*$-James space.

\proclaim Proposition 8.  ${T{}_J^*}$ has a monotone
Schauder basis.

\noindent{\bf Proof}.  Let $(e_j)_{j=1}^\infty$ denote the
canonical unit vector basis for $c_{00}\subset {T{}_J^*}$.

By (3) it is clear that
$$\left\|\sum_{j=1}^na_je_j\right\|_{T{}_J^*}\le
\left\|\sum_{j=1}^{n+1}
a_je_j\right\|_{T{}_J^*}$$ for any sequence of scalars
$(a_j)_{j=1}^{n+1}$ and hence, by [27, proposition 1.a.3],
the sequence $(e_j)_{j=1}^\infty$ is a monotone Schauder
basis for $T{}_J^*$.

Our next proposition shows that the norms on $T^*$ and
$T{}_J^*$ behave in the same way with respect to normalised
block basic sequences.

\proclaim Proposition 9.  Let $(u_n)_{n=1}^\infty$ denote a
normalised block basic sequence in $T{}_J^*$.  For any
sequence of scalars $(a_j)_{j=1}^\infty$ and for any
positive integer $n$ we have
$$\left\|\sum_{j=1}^na_ju_j\right\|_{T{}_J^*}\le
\left\|\sum_{j=1}^n\bigl(|a_j|+|a_{j+1}|\bigr)t_j
\right\|_{T^*}$$

\noindent{\bf Proof}.  We first fix $n$ and let
$$\sum_{l=1}^\infty b_le_l=
\sumj 1 n a_ju_j\quad\hbox{where }
u_j=\sumi {{k_j+1}} {{k_{j+1}}} a_{ij}e_i$$
and $(k_j)_{j=1}^\infty$ is a strictly increasing sequence
of integers with $k_1=0$.  Let $1\le p_1<p_2<\cdots<p_{2k}\le
k_{s+1}$ denote an increasing sequence of
positive integers.

By (2)
$$\eqalign{
\left\|\sumj 1 k(b_{p_{2j-1}}-b_{p_{2j}})t_j\right\|_{T^*}
&\le\left\|\Biggl\|\sum_{\scriptstyle j\atop
\scriptstyle  p_{2j-1}\le k_{2}}
(b_{p_{2j-1}}-b_{p_{2j}})t_j\right.\Biggr\|_{T^*}t_1\cr
&\left.+\sum_{l=2}^s\Biggl\|\sum_{\scriptstyle j\atop
\scriptstyle  k_l<p_{2j-1}\le k_{l+1}}
(b_{p_{2j-1}}-b_{p_{2j}})t_j\Biggr\|_{T^*}t_l
\right\|_{T^*}\cr}$$
Let $a_0=a_{s+1}=0$. Then
$$\eqalign{
\left\|\sumj 1
k(b_{p_{2j-1}}-b_{p_{2j}})t_j\right\|_{T^*}
&\le \left\|\suml 1
{{s}}\bigl(\|a_lu_l\|_{T{}_J^*}
+|a_{l+1}|\bigr)t_l\right\|_{T^*}\cr
&\le\left\|\suml 1 {{s}}
\bigl(|a_l|+|a_{l+1}|\bigr)t_l\right\|_{T^*}\cr}$$
Hence
$$\sup_{\scriptstyle p_1<\cdots<p_{2k}\atop
\scriptstyle  k}\left\|
\sum_{j=1}^k(b_{p_{2j-1}}-b_{p_{2j}})t_j\right\|_{T^*}$$
$$=\left\|\sumj 1 n a_ju_j\right\|_{T{}_J^*}\le
\left\|\sumj 1 n \bigl(|a_j|+|a_{j+1}|\bigr)t_j
\right\|_{T^*}$$

\proclaim Corollary 10.  If $(u_j)_{j=1}^\infty$ is a
normalised block basic sequence in $T{}_J^*$ then
$$\left\|\sumj n {{2n}}
a_ju_j\right\|_{T{}_J^*}\le4\sup_{n\le j\le 2n}|a_j|$$

\noindent{\bf Proof}.  It suffices to apply (1) and
proposition~9.
\bigskip
A sequence of vectors $(u_j)_j$ in a Banach
space is said to have a lower $q$-estimate if there exists
$c>0$ such that
$$\left\|\sumj n m a_iu_j\right\|
\ge c\left(\sumj n m |a_j|^q\right)^{1/q}$$
for any positive integers $n$ and $m$.
\bigskip
\proclaim Corollary 11.  No normalised block basic sequence
in $T{}_J^*$ satisfies a lower $q$ estimate for any
$q<\infty$.

\noindent{\bf Proof}.  Otherwise, by
corollary~10, there would exist $q<\infty$ and $c>0$ such
that
$$\left(\sumj n {{2n}} |a_j|^q\right)^{1/ q}\le c \sup_{n\le
j\le 2n}|a_j|$$
for all sequences of scalars $(a_j)_{j=1}^\infty$.

If $n=2^{j}$, $m=2^{j+1}$ and $a_j=1$ for all $j$ this
would imply
$$2^{j/q}\le c$$ for all $j$ and this is impossible.

\proclaim Corollary 12.  The sequence
$(e_j)_{j=1}^\infty$ is
a shrinking basis for $T{}_J^*$.

\noindent{\bf Proof}.  Otherwise there would exist
$\phi\in(T{}_J^{*})'$, $(u_j)_{j=1}^\infty$ a normalised
block basic sequence in $T{}_J^*$, and $\delta>0$ such that
$\phi(u_j)\ge\delta$ for all $j$.

Let $\dis\alpha_j={1\o2^n}$ for $2^n<j\le2^{n+1}$,
$n=1,\ldots$ We have $\sumj 1 \infty \alpha_j=\infty$ and
$\alpha_j>0$ for all $j$.  On the other hand, by
corollary~10,
$$\left\|\sumj {{2^n+1}}
{{2^{n+1}}}\alpha_ju_j\right\|_{T{}_J^*}
\le4\sup_{2^n<j\le2^{n+1}}|\alpha_j|\le{4\o2^n}$$
and
$$\eqalign{
\sumj 1 \infty\left\|\sumj {{2^n+1}}
{{2^{n+1}}}\alpha_ju_j\right\|_{T{}_J^*}
&\le4\sumn 1 \infty\sup_{2^n<j\le2^{{n+1}}}|\alpha_j|\cr
&\le4\sumn 1 \infty{1\o2^n}<\infty\cr}$$
Hence $\sumj 1 \infty\alpha_ju_j\in T{}_J^*$.

However,
$$\lim_{n\to\infty}\phi\left(\sumj 1
n\alpha_ju_j\right)
\ge\delta\lim_{n\to\infty}\sumj 1
n\alpha_j=\infty.$$
This is impossible and shows that the basis is shrinking.
\bigskip
By corollary~12 and [27, proposition 1.b.1 and 1.b.2] the
biorthogonal functions $(e_j^*)_{j=1}^\infty$ form a
Schauder basis for $(T{}_J^*)'$ and, moreover,
$(T{}_J^*)''$ can
be identified with the space of all sequences
$({a_j})_j$ such
that
$$\sup_n\left\|\sumj 1 na_je_j\right\|_{T{}_J^*}<\infty.$$
The correspondence is given by
$$x^{**}\in(T{}_J^*)''
\longleftrightarrow\bigl(x^{**}(e_j^*)\bigr)_{j=1}^\infty$$
and we have
$$\|x^{**}\|=\sup_n\left\|\sumj 1
nx^{**}(e_j^*)e_j\right\|_{T{}_J^*}.$$

\proclaim Proposition 13.  $T{}_J^*$ is not reflexive.

\noindent {\bf Proof}.  Let $w_n=\sumj 1 n e_j\in T{}_J^*$.
Let $b_j=1$ for $j\le n$ and $b_j=0$ for $j>n$.  Since
$$\sum_{\scriptstyle j\atop\scriptstyle
p_1<p_2<\cdots<p_{2k}}(b_{p_{2j-1}}-b_{2j})e_j
=\cases{
0&if $p_1>n$ or $p_{2k}\le n$\cr
e_j&if $p_{2j-1}\le n\le p_{2j}$ for some $j$,\cr
0&if $p_{2j}\le n<p_{2j+1}$ for some $j$,\cr}$$
it follows that $\|w_n\|_{T{}_J^*}=1$ for all $n$.
If $T{}_J^*$ were reflexive then the sequence
$\{w_n\}_{n=1}^\infty$ would contain a
subsequence which was
weakly convergent to some $w\in T{}_j^*$.
Since $e_m^*(w_n)=1$
for all $n\ge m$ it follows that $e_m^*(w)=1$ for all $m$.
If $\dis w=\summ 1 \infty \beta_me_m\in T{}_J^*$ then
$\|\beta_me_m\|_{T{}_J^*}\to0$ and hence $|\beta_m|\to0$ as
$m\to\infty$.  Since $\beta_m=1$ for all $m$ this is
impossible and completes the proof.
\bigskip
We now describe $(T{}_J^*)''$.
If  $(a_j)_{j=1}^\infty$ is a sequence of scalars and
$\dis\sup_n\left\|\sumj 1 n a_j e_j\right\|_{T{}_J^*}
=M<\infty$
then we claim that $\dis\lim_{j\to\infty}a_j$ exists.
Otherwise, there exists $\delta>0$ and a strictly
increasing sequence of
positive integers $(p_j)_j$ such that
$|a_{p_{2j-1}}-a_{p_{2j}}|\ge\delta>0$ for all $j$.  Hence
$$\sup_k\left\|\sumj 1
k(a_{p_{2j-1}}-a_{p_{2j}})t_j\right\|_{T^*}\le M.$$
Since the basis in $T^*$ is 1-unconditional we have
$$\sup_k\left\|\sumj 1 k\delta t_j\right\|\le
\sup_k
\left\|\sumj 1
k|a_{p_{2j-1}}-a_{p_{2j}}|t_j\right\|\le M.$$
Since $T^*$ is reflexive this implies that
the sequence $\dis\left\{\sumj 1 k
t_j\right\}_{k=1}^\infty$ has a weak Cauchy subsequence and
the proof of the previous proposition can now be
adapted to show that this is impossible.  Hence we have
established our claim.

In proving proposition 13 we showed that
$\dis\sup_n\left\|\sumj 1 n e_j\right\|_{T{}_J^*}<\infty$.
Let $x_0^{**}\in(T{}^*_J)''$ be given by $x_0^{**}(e_n^*)=1$
for all $n$.  If $x^{**}$ is an arbitrary vector in
$(T{}_J^*)''$ then
$\dis\lim_{j\in\infty}x^{**}(e_j^*)=\alpha(x^{**})$ exists.
Let $y^{**}=x^{{**}}-\alpha(x^{**})x_0^{**}$. It follows
that $y^{**}\in T{}_J^*$ and $(T{}_J^*)''\cong
T{}_J^*\oplus{\bf C}x_0^{**}$. Hence $(T{}_J^*)''$ is
separable and by [27, theorem 1.c.12] the basis
$(e_j)_{j=1}^\infty$ in $T{}_J^*$ is not unconditional.

The above also shows that $T{}_J^*$ and all its duals are
quasi-reflexive and hence all have the Radon-Nikodym
property [16, p. 219].  Moreover, $T{}_J^*$ and all its
higher duals have a basis and hence the approximation
property.

The following proposition may also be proved by using the
method used for Tsirelson's space in [3] and corollary~10.

\proclaim Proposition 14.  Continuous polynomials on
$T{}_J^*$ are weakly continuous on bounded sets.

\noindent{\bf Proof}.  Since $T{}_J^*$ and $(T{}_J^*)''$ are
both separable it follows by [27, Theorem 2.e.7] that
$l_1\not\hookrightarrow T{}_J^*$.  Hence by proposition~2,
it suffices to show that each continuous polynomial on
$T{}_J^*$ is weakly sequentially continuous at the origin.
By proposition~3, it suffices to show that no spreading
model built on a normalised weakly null sequence in
$T{}_J^*$ has a lower $q$ estimate for any $q<\infty$.

Suppose $(u_j)_j$ is a normalised weakly null sequence in
$T{}_J^*$ which has a spreading model having a lower $q$
estimate for some $q<\infty$.  This means that there exists
a Banach space $X$ with an unconditional basis, $(f_j)_j$,
{\it the spreading model}, and a subsequence of $(u_j)_j$,
$(u_{n_j})_{j=1}^\infty$, such that for all $\epsilon>0$ and
all $k$ there exists $N=N(\epsilon,k)$ such that for all
$N<n_1<n_2<\cdots<n_k$ and for all scalars with
$\dis\sup_i|a_i|\le1$ we have
$$\left|\Biggl\|\sumj 1 k a_jf_j\Biggr\|_X
-\Biggl\|\sum_{j=1}^k a_j
u_{n_j}\Biggr\|_{T{}_J^*}\right|<\epsilon$$
(see for instance [10, 18, 19, 25])
and that there exists $c>0$ such that
$$\left\|\sumj 1 k a_jf_j\right\|_X\ge c
\left(\sumj 1 k |a_j|^q\right)^{1/q}$$
for any sequence of scalars $(a_j)_j$ and for all $k$.

By choosing, if necessary, a further subsequence of
$(n_j)_j$ we may suppose that
$$u_{n_j}=x_{n_j}+y_{n_j}$$
where $(x_{n_j})_j$ is a block basic sequence in $T{}_J^*$
and $\dis\sumj 1 \infty\|y_{n_j}\|_{T{}_J^*}\le{1\o4}$.

Since $\|u_{n_j}\|=1$ for all $j$ this implies that
$\dis\|x_{n_j}\|\ge{3\o4}$ for all $j$.
By corollary 10 we have for $k$
sufficiently large and $|a_j|\le1$.
$$\eqalign{
4\sup_{k+1\le j\le 2k}|a_j|
&\ge\left\|\sumj {{k+1}} {{2k}}
a_jx_{n_j}\right\|_{T{}_J^*}\cr
&\ge\left\|\sumj {{k+1}} {{2k}}
a_ju_{n_j}\right\|_{T{}_J^*}
-\left\|\sumj {{k+1}} {{2k}}
a_jy_{n_j}\right\|_{T{}_J^*}\cr
&\ge\left\|\sumj {{k+1}} {{2k}}
a_jf_{j}\right\|_X
-\sup_{1\le j\le k}|a_j|
\sumj {{k+1}} {{2k}}
\|y_{n_j}\|_{T{}_J^*}-\epsilon\cr
&\ge c\left(\sumj {k+1} {2k}|a_j|^q\right)^{1/q}-{1\o4}
\sup_{k+1\le j\le2k}|a_j|-\epsilon.\cr}$$
Since $\epsilon>0$ was arbitrary this implies that for any
$k$ and any sequence $(a_j)_j$ with $|a_j|\le1$ we have
$$c\left(\sum_{j=1}^k|a_j|^q\right)^{1/q}\le 5\sup_{1\le
j\le k}|a_j|.$$
Letting $a_j=1$ for all $j$ this implies $\dis ck^{1/ q}\le
5$
for all positive integers $k$.  This contradiction
shows that
no normalised weakly null sequence in $T{}_J^*$ has a
spreading model with a lower $q$ estimate for some
$q<\infty$ and completes the proof.

\proclaim Proposition 15.  $T{}_J^*$ is a $Q$-reflexive
Banach space and ${\cal P}(^nT{}_J^*)$ has RNP for all $n$.

\noindent{\bf Proof}.  Since $T{}_J^*$ is a quasi-reflexive
space with basis, it follows that $(T{}_J^*)'$ and
$(T{}_J^*)''$ have the approximation property and RNP.  The
remaining hypothesis required for the application of
Theorem~5 is satisfied by proposition~14.

\bigskip
The hypotheses of theorem 5 are also satisfied by any
subspace of $(T{}_J^*)^{(2n)}$, the
$2n^{\hbox{\sevenrm th}}$ dual of
$T{}_J^*$, which has the approximation property.

The mapping
$$U:(T{}_J^*)''\to T{}_J^*$$
$$U(x^{**})=(-\lambda,x^{**}(e^*_1)-
\lambda,x^{**}(e^*_2)-\lambda,\ldots)$$
where $\lambda=\alpha(x^{**})$ is linear isomorphism from
$(T{}_J^*)''$ onto $T{}_J^*$.

Consequently, since $E$ is $Q$-reflexive, ${\cal
P}(^nT{}_J^*)$ and ${\cal P}(^nT{}_J^*)''$ are isomorphic
for all $n$, (but not under the canonical mapping since
$T{}_J^*$ is not reflexive).
This property, which is shared by many, but possibly not all,
quasi-reflexive spaces, suggested the terminology
$Q$-reflexive spaces.

If $E$ is a Banach space we let ${\cal H}_b(E)$ denote the
space of {\bf C}-valued holomorphic  functions on $E$ which
are bounded on bounded sets and endowed with the topology
$\tau_b$ of uniform convergence on bounded sets.

\proclaim Proposition 16.  $\dis\bigl({\cal
H}_b(T{}_J^*),\tau_b\bigr)''\cong \Bigl({\cal
H}_b\bigl((T{}_J^*)''\bigr),\tau_b\Bigr).$

\noindent{\bf Proof}.  It suffices to apply theorem 12 of
[28] and proposition 15.
\bigskip

Further applications to spaces of holomorphic functions are
also possible and we will discuss these in a further paper.

Finally, we have recently obtained preprints by J. A. Jaramillo,
A. Prieto, I. Zalduendo [24] and by M. Valdivia [30], in
which the bidual of $\pne$ is discussed.

\vglue.5truein
\centerline{\bf References}
\bigskip
\item{\bf[1]} R. Alencar, {\sl On reflexivity and basis for
$\pme$}, P.R.I.A., {\bf85A}, 2, 1985, 131--138.
\item{\bf[2]} R. Alencar, {\sl An application of Singer's
theorem to homogeneous polynomials}, AMS, Contemporary
Math.,  Ed. Bor-Luh Lin, W. B. Johnson, {\bf 144}, (1993),
1--8.
\item{\bf[3]}  R. Alencar, R. M. Aron, S. Dineen. {\sl A
reflexive space of holomorphic functions in infinitely many
variables,} P.A.M.S.,{\bf90}, 1984, 407--411.
\item{\bf[4]}R. Alencar, R. M. Aron, G. Fricke, {\sl Tensor
products of Tsirelson's Space}, Illinois J. Math.,  {\bf31},
1, 1987, 17--23.
\item{\bf[5]} R. M. Aron, {\sl Compact polynomials and
compact differentiable mappings between Banach spaces},
Seminaire P. Lelong 1974--1975, Springer Verlag Lecture
notes, {\bf 524}, (1976), 213--222.
\item{\bf[6]} R. M. Aron \& P. Berner, {\sl A Hahn Banach
extension theorem for analytic mappings}, Bull. Soc. Math. France,
{\bf 106}, 1978, 3--24.
\item{\bf[7]} R. M. Aron, J. Diestel, A. K. Rajappa, {\sl A
characterization of Banach spaces containing $l_1$, Banach
Spaces}, Proceedings of the Missouri Conference, June 1984,
Springer Verlag lecture Notes in Math., {\bf1166}, 1985,
1--3.
\item{\bf[8]} R. M. Aron, J. Globevnik, {\sl Analytic
functions on $c_0$}, Revista Matem\'atica, {\bf2}, 1989,
27--33.
\item{\bf[9]} R. M. Aron, C. Herves, M. Valdivia, {\sl
Weakly continuous mappings on Banach spaces}, J. F.
Analysis, {\bf52}, 2, 1983, 189--204.
\item{\bf[10]}  B. Beausamy, J.-T.  Laprest\'e, {\sl
Mod\`eles \'etal\'es des espaces de Banach,} Hermann, Paris,
1984.
\item{\bf[11]} P. G. Casazza, P. G. Lin, R. H. Lohman,
{\sl On non-reflexive spaces which contain no $c_0$ or
$l_p$},  Can. J. Math., {\bf32},  1980, 1382--1389.
\item{\bf[12]} P. G. Casazza, T. J. Shura, {\sl Tsirelson's
spaces}, Springer-Verlag Lecture Notes in Math.
{\bf1363}, 1989.
\item{\bf[13]}  A. M. Davie, {\sl Quotient algebras of
uniform algebras,}  J. London Math. Soc. {\bf7}, 1973,
31--40.
\item{\bf[14]} A. Defant, {\sl A duality theory for locally
convex tensor products}, Math. Z., {\bf190},  1985, 45--53.
\item{\bf[15]} J. Diestel, {\sl Sequences and series in
Banach spaces}, Springer-Verlag, Graduate Texts in
Mathematics, 1984.
\item{\bf[16]}  J. Diestel and J. Uhl, {\sl Vector
Measures},  American Math Surveys, 1977.
\item{\bf[17]} S. Dineen, {\sl Complex Analysis in Locally
Convex Spaces}, North Holland Math.~Studies,  {\bf59}, 1981.
\item{\bf[18]} J. D. Farmer, {\sl Polynomial reflexivity in
Banach spaces}, preprint, 1992.
\item{\bf[19]} J. Farmer, W. B. Johnson, {\sl Polynomial
Schur and polynomial Dunford-Pettis properties},
AMS, Contemporary
Math.,  Ed. Bor-Luh Lin, W. B. Johnson,
{\bf 144}, (1993),
95--105.
\item{\bf[20]}  J. E. Gilbert, T. J. Leih, {\sl
Factorization, tensor products and bilinear forms in Banach
space theory,}  notes in Banach spaces, Ed. H.~Elton~Lacey,
University of Texas Press, Austin and London, 1980,
182--305.
\item{\bf[21]} M. Gonzalez, J. M. Guti\'errez, {\sl
Unconditionally converging polynomials on Banach spaces},
preprint 1992.
\item{\bf[22]}  A. Grothendieck, {\sl Produits tensoriels
topologiques et espaces nucl\'eaires}, Memoirs, A.M.S.,
{\bf16}, (1955), (reprinted 1966).
\item{\bf [23]}  J. Gutierrez, {\sl Weakly continuous functions
on Banach spaces not containing $\ell_{1}$}, Proc. A.M.S.,
{\bf 119}, no. 1 (1993), 147--152.
\item{\bf [24]}  J.A. Jaramillo, A. Prieto, and I. Zalduendo,
{\sl The bidual of a space of polynomials and polynomials on the
bidual of a Banach space}, preprint.
\item{\bf[25]} H. Knaust, E. Odell, {\sl Weakly null
sequences with upper $l_p$-estimates},
Ed. E. Odell, H. Rosenthal,
Functional Analysis
Proceedings,
Springer-Verlag Lecture Notes in Mathematics, {\bf1470},
1991, 85--107.
\item{\bf[26]}  R. H. Lohman, P.G. Casazza, {\sl A general
construction of spaces of the type of R. C. James},
Can.~J.~Math.,~{\bf27}, 1975, 1263--1270.
\item{\bf[27]} J. Lindenstrauss, L. Tzafriri, {\sl
Classical
Banach spaces~I, Sequences Spaces}, Springer-Varlag,  1977.
\item{\bf[28]} A. Prieto, {\sl A bidual of spaces of
holomorphic functions in infinitely many variables},
P.R.I.A., {\bf92A}, 1, 1992, 1--8.
\item{\bf[29]}  B. Sims, D. Yost, {\sl Banach spaces with many
projections}, in ``Miniconference on operator theory and partial
diffential equations,'' Proc. Centre Math. Anal. Austral. Nat. Univ. 
{\bf 14} (1986), 335-342.
\item {\bf[30]}  M. Valdivia,  {\sl  Banach spaces of
polynomials without copies of $\ell^1$}, preprint.
\item{\bf[31]} I. Zalduendo, {\sl An estimate for
multilinear forms on $l^p$ spaces}, Proc. R. I. Acad., 
Sect. A {\bf 93}, no. 1 (1993), 137-142.

\line{Department of Mathematics\hfill Department of
Mathematics}
\line{University College Dublin\hfill Kent State
University}
\line{Belfield\hfill Ohio 44242}
\line{Dublin 4\hfill U.S.A.}
\line{Ireland\hfill }
\medskip
\line{e-mail: sdineen@irlearn.bitnet\hfill e-mail
aron@mcs.kent.edu }

\bye